\documentclass{amsart}
\usepackage{amsmath, amsthm}

\title{Fixed Point Theorem for Discontinuous mappings on PN Spaces}

\author{Mohd Rafi}
\address{Faculty of Engineering and Computer Science, The University of Nottingham Malaysia Campus, Jalan Broga, 43500 Semenyih, Selagor Darul Ehsan, Malaysia}
\email{mdrafzi@yahoo.com}


\newtheoremstyle{theorem}%name
  {10pt}		  % space above
  {10pt}  % space below
  {\sl}  % bofy font
  {\parindent}     % ident - empty=no indent,  \parindent= paragraph indent
  {\bf}  % thm head font
  {. }    % punctuation after thm head
  { }    % space after thm head: `` ``=normal \newline=linebreak
  {}     % thm head specification
\theoremstyle{theorem}

\newtheoremstyle{defi}%name
  {10pt}		  % space above
  {10pt}  % space below
  {\rm}  % bofy font
  {\parindent}     % ident - empty=no indent,  \parindent= paragraph indent
  {\bf}  % thm head font
  {. }    % punctuation after thm head
  { }    % space after thm head: `` ``=normal \newline=linebreak
  {}     % thm head specification
\theoremstyle{defi}

\begin{document}

\maketitle

\begin{abstract}
We present a study on strong t-continuity and measure of discontinuity on PN spaces. As an application, we prove a fixed point theorem for a discontinuous self mapping on PN spaces by means of measure of discontinuity.
\end{abstract}

\section{Introduction}
The notion of probabilistic normed space (briefly, a PN space) was introduced by Serstnev. The theory of PN spaces has been dormant after its initial application untill it was recently put on new basis by Alsina, Schweizer and Sklar \cite {Alsina}. In 1997, B. L. Guillen et al further studied some classes of probabilistic normed spaces. In this paper, the author has considered Serstnev spaces and obtained some results using minimum t-norms. The results in this paper would be much improved if they were proved for arbitrary t-norms. Recently results about bounded linear operators and continuous operators in PN spaces have been studied using arbitrary t-norms \cite{Iqbal}.

A study on fixed points of contraction mappings on probabilistic metric spaces (PM spaces) was initiated by V. M. Sehgal and A. T. Bharucha-Reid \cite{Sehgal} and later developed by numerous authors, as can be realized upon consulting the list of references given at \cite{Cons} and \cite{Hadzicpap}. One the other hand, O. Hadzic \cite{Hadzic} initiated fixed point theory in random normed spaces. Recently, in \cite{Rafi}, the authors have studied approximate fixed point theorems on PM and PN spaces respectively which generalizes the results obtained in \cite{Branzei} and \cite{Tijs}. 

The purpose of this paper is to prove the existence of fixed points for discontinuous self mappings in PN spaces, a study which carried out by L. J. Cromme and I. Diener \cite{Cromdien} on normed spaces.

In the sequel we generally follow the notation and terminology of \cite{Schweizer}.

\section{Preliminaries} 

A {\em distance distribition function} (briefly, a d.d.f.) is a function $F$ from the extended positve real line $\overline{\mathbb{R}}=[0,+\infty]$ into the unit interval $I=[0,1]$ that is nondecreasing and satisfies $F(0)=0$, $F(+\infty)=1$. We normalize all d.d.f.'s to be left-continuous on the unextended positive real line $\mathbb{R}= ]0,+\infty[$. Particularly, for every $a$ in $[0, +\infty[$, the functions $\varepsilon_a$ is d.d.f. defined by

\[  \varepsilon_a(x)= \left\{\begin{array}{c l}
    0, & x\leq a;\\
    1, & x>a.
  \end{array}
\right. \]

The set of all d.d.f.'s will be denoted by $\Delta^+$. The set $\Delta^+$ is partially ordered by the ususal pointwise partial ordering of functions, i.e., $F\leq G$ whenever $F(x)\leq G(x)$, for all $x$ in $\mathbb{R}$. The maximal element for $\Delta^+$ in this order is the d.d.f. $\varepsilon_0$.

A sequence $\{F_n\}$ in $\Delta^+$ is said to be weakly convergent to $F \in \Delta^+$ (shortly, $F_n \stackrel{w}\to F$) if $\lim_{n\to\infty}F_{n}(x)=F(x)$ for every continuity point $x$ of $F$.

A {\em triangle function} is a binary operation on $\Delta^+$ 
that is commutative,associative, nondecreasing in each place and has
$\varepsilon_0$ as identity. Continuity of a triangle function means
uniform continuity with respect to the natural product topology on $\Delta^+\times\Delta^+$.

The typical triangle function is the
operation  $\tau_T$ which given by
$$\tau_T(F,G)(x)=\sup_{u+v=x}T(F(u),G(v))$$
for all $F,G\in\Delta^+$ and all $x>0$
(\cite{Schweizer}, section 7.2 and 7.3). Here  $T$ is a $t$-{\em norm},
i.e., $T$ is binary operation on $[0,1]$  that is commutative,
associative, nondecreasing in each place and has 1 as identity.

The most important t-norms are the function $W$, $Prod$, and $M$
which are defined, respectively, by

\begin{eqnarray*}
W(a,b)&=&\max\{a+b-1,0\},\\
Prod(a,b)&=&ab,\\
M(a,b)&=&\min(a,b).
\end{eqnarray*}

\newtheorem{defn}{Definition}[section]
\begin{defn}(\cite{Alsina})
A PN space is a quadruple $(V,\nu,\tau, \tau^*)$ where $V$ is a real linear vector space, $\tau$ and $\tau^*$ are continuous triangle functions and the probabilistic norm $\nu$ is a mapping from $V$ into $\Delta^+$ such that, for all $p,q$ in $V$, the following conditions hold (we use $\nu_{p}$ instead of $\nu(p)$):\\
(N1) $\nu_{p}=\varepsilon_0$, if and only if $p=\vartheta$ (where $\vartheta$ is a null vector in $V$),\\
(N2) $\nu_{-p}=\nu_p$\\
(N3) $\nu_{p+q}\geq \tau(\nu_p,\nu_q)$,\\
(N4) $\nu_{p} \leq \tau^*(\nu_{\lambda p},\nu_{(1-\lambda)p})$ for every $\lambda \in [0,1]$.
\end{defn}

Since $\tau$ is continuous, the system of strong neighborhoods $\{N_{p}(t)|p\in V, t>0\}$, where
$$N_{p}(t)=\{q\in V|v_{p-q}(t)>1-t\},$$
determines a topology on $V$, called the {\em strong topology}.

\begin{defn}(\cite{Guillen})
Let $(V,\nu,\tau, \tau^*)$ be a PN space and $A$ is a nonempty subset of
$V$. The probabilistic diameter of $A$, $R_A$ is defined as
\[ R_{A}(x)=\left\{\begin{array}{cl}
l^{-}\varphi_{A}(x),&x\in[0,+\infty);\\
1,&x=+\infty,
\end{array}
\right. \]
where $l^{-}f(x)$ denotes the left limit of the function $f$ at $x$ and \[ \varphi_{A}(x)=\inf_{p\in A}\nu_{p}(x)\].
\end{defn}

As a consequence of Definition 2.2, we have $\nu_{p}\geq R_{A}$.

\section{MAIN RESULTS}

Now we propose the following {\em strong t-continuity} in PN spaces.

\begin{defn}
In a PN space $(V,\nu,\tau, \tau^*)$, a map $f\colon V\to V$ is said to be strong $t$-continuous, $(t>0)$ provided for each $p\in V$ admits a strong $t'$-neighborhood $N_{p}(t')$, $t'>0$ such that

$$R_{f(N_{p}(t'))}(t)>1-t.$$
\end{defn}

The following theorem is a direct consequence of the Definition 3.1 above.

\newtheorem{theo}{Theorem}[section]
\begin{theo}
 Let $(V,\nu,\tau, \tau^*)$  be a PN space such that $\tau=\tau_{M}$ and $A\subseteq V$. Let $f\colon A \to A$ be a strong $t$-continuous mapping where $t>0$, then for every distinct $p,q\in A$,
$$\nu_{f(p)-f(q)}(t)>1-t.$$ 
\end{theo}
{\bf Proof}\quad Since $f$ is a strong $t$-continuous mapping, for $p,q\in A$, there exist strong neighborhoods $X=N_{p}(t'),Y=N_{q}(t")$ for every $t',t">0$ such that 

$$\nu_{f(p)}(t)\geq R_{f(X)}(t)>1-t$$
and
 
$$\nu_{f(q)}(t)\geq R_{f(Y)}(t)>1-t$$.

It follows that,
\begin{eqnarray*} 
\nu_{f(p)-f(q)}(t)&\geq& \tau_{M}(\nu_{f(p)},\nu_{f(q)})(t)\\
			&=&\sup_{s\in (0,1)}M(\nu_{f(p)}(st),\nu_{f(q)}(t-st))\\
			&\geq& M(\nu_{f(p)}(t),\nu_{f(q)}(t))\\
			&>&Min(1-t,1-t)=1-t.
\end{eqnarray*}

The existence and uniqueness of fixed points for classes of mappings on PN spaces depending on some kind of continuity assumption. Without this assumption, the theorem could be not valid; it needs not be completely false, though, see \cite{Cromdien} and \cite{Crom}.

The discontinuous of a self mapping can be interpreted several ways. Here, we consider the jump discontinuity, i.e.,$\lim_{x \to a}f(x)\neq f(a)$. For this purpose, we propose the following measure of discontinuity in PN spaces (a generalization of the measure of discontinuity in normed spaces, see \cite{Crom}). 

\begin{defn}  
Let $(V,\nu,\tau, \tau^*)$   be a PN space, $A\subseteq V$ be precompact and $f\colon A \to A$ any map. We define the following {\em probabilistic measure of discontinuity}:
$$\psi_{f}(t)=\inf_{p\in A}\liminf_{\delta \to 0}\inf_{q\in N_{p}(\delta)}\nu_{f(p)-f(q)}(t), t>0$$.
\end{defn}

It can be noted that if  $f$ is strong $t$-continuous self mapping then $\psi_{f}=\varepsilon_0$..

Any (discontinuous) function $f\colon A\to A$, can be transformed into a continuous point-to-set mapping $T_{f}$ as follows: for $p\in A$, we set
$$T_{f}(p)=\left\{f(p_i)\colon \lim_{i\to\infty}\nu_{p_{i}-p}=\varepsilon_0\right\}$$,

i.e., the set $T_{f}(p)$ assigned to a point $p\in A$ includes all accumulation points of sequences $\left\{f(p_{i})\colon \lim_{i\to \infty}\nu_{p_{i}-p}=\varepsilon_0\right\}$. Whence the point-to-set mapping \( p \stackrel{T}{\rightarrow} conT_{f}(p)\)(where $conT_{f}(p)$ is convex hull of set $T_{f}(p)$) is upper semi-continuous and may be interpreted as 'continuoufication' or smooting of the discontinuous original function $f$. 

\newtheorem{lemma}{Lemma}[section]
\begin{lemma}

 Let $(V,\nu,\tau, \tau^*)$ be a PN space, $A\subseteq V$  be precompact and $f\colon A \to A$  any map. Let $q_{i}\to p^*$ be a convergent sequence in $A$ and suppose for every $i$, there exists a convergent sequence $p_{ij}\to p^{*}_{i}$ in $A$ with
$$\lim_{j\to \infty}\nu_{fp_{ij}-q_i}=\varepsilon_0$$. 

Suppose further that $p^{*}_{i}\to p^*$  for some $p^{*}\in V$. Then there exists a sequence $r_{j}\to p^*$ in A such that
$$\lim_{j\to \infty}\nu_{fr_{j}-q^*}=\varepsilon_0$$
\end{lemma}
{\bf Proof}\quad By hypothesis, we have

$$\lim_{i\to \infty}\nu_{q_{i}-q^*}=\varepsilon_0, \lim_{i\to \infty}\nu_{p^{*}_{i}-p^{*}}=\varepsilon_0$$ and $$\lim_{i\to \infty}\nu_{p_{ij}-p^*}=\varepsilon_0$$.
Let $r_{k}=p_{i_{k}j_{k}}$. By (PN2) and (PN3) we have,

\begin{eqnarray*}
\nu_{r_{k}-p^*}&\geq&\tau(\nu_{r_{k}-p^{*}_{i_k}},\nu_{p^{*}_{i_k}-p^*})\\
		    &\geq&\tau(\varepsilon_0,\varepsilon_0)=\varepsilon_0\\
and\\
\nu_{fr_{k}-q^*}&\geq&\tau(\nu_{fr_{k}-q_{i_k}}.\nu_{q_{i_k}-q^*})\\
		     &\geq&\tau(\varepsilon_0,\varepsilon_0)=\varepsilon_0
\end{eqnarray*}
as $k\to \infty$. Thus $r_{k}\to p^*$ and $\lim_{k\to \infty}\nu_{fr_{k}-q^*}=\varepsilon_0$, as required.

\begin{lemma}

 Let $(V,\nu,\tau, \tau^*)$ be a PN space, $A$ be a precompact and convex subset of $V$ and let $f \colon A\to A$ be any map. We define a set valued map $T_{f}(p)$ on $A$:

$$T_{f}(p)=\{q\in A|\exists p_{i}\to p, p_{i}\neq p, \lim_{i\to\infty}\nu_{fp_{i}-q}=\varepsilon_0\}$$

If $V=\{p\}$ we set $T_{f}(p)=\{p\}$. Then there exists a point $p^{*}\in A$ such that $p^{*}\in conT_{f}(p^*)$.
\end{lemma}
{\bf Proof}\quad Let $F(p)=conT_{f}(p)$. We show that the graph $G=\{(p,q)|p\in A, q\in F(p)\}$ is closed subset of $A\times A$. Let $(p_{i},q_{i})\to (p,q)$ be a convergent sequence in $A\times A$ with $q_{i}\in F(p_i)$. If $p_i$ equals $p$ infinitely often then it is trivial. Without loss of generality we assume $p_{i}\neq p,\quad\forall i=1,2,\ldots,n$. By theorem of Caratheodory, since $q_{i}\in conT_{f}(p_i)$, we can for every $q_i$ choose $n+1$ points $p_{i_0},\ldots,p_{i_n}\in T_{f}(p_i)$ and scalars $\lambda_{i_0},\ldots, \lambda_{i_j}\geq 0$ with $\sum\lambda_{i_j}=1$, such that
$$q_{i}=\sum_{i=1}^{\infty} \lambda_{i_j}p_{i_j}$$.

Since $\lambda_{i}=(\lambda_{i_0},\ldots, \lambda_{i_n})$ is contained in the compact standard simplex $S_n$ and the $p_{i_j}$ are contained in the precompact set $A$ we can consider a subsequence $\{q_{i_j}\}\subset\{q_i\}$ with $\nu_{q_{i_j}-q}\stackrel{w}\to\varepsilon_0$ such that the corresponding subsequences $\{\lambda_{k_i}\}_{i=1}^{\infty},\quad\j=1,\ldots,n$ and $\{p_{k_i}\}_{i=1}^{\infty},j=1,\ldots,n$, converges to points $\lambda_k$ and $p_k$ respectively. It follows that $q$ can be written as
$q=\sum_{k=0}^{n}\lambda_{k}p_k$.
Thus, $q\in con\{p_0,\ldots,p_n\}$. It remains to show that $p_k\in T_{f}(p)$ for $k=0,\ldots,n$. For fixed $k$ we set $s_{i}=p_{k_i}$ and $s^{*}=p_k$. Since $s_{i}\in T_{f}(p_i)$ there exists a sequence $t_{ij}\to p_{i}$, $t_{ij}\neq p_i$ such that $\nu_{f(t_ij)-s_i}\stackrel{w}\to\varepsilon_0$. Furthermore, $p_i\to p$. Thus from Lemma 3.1, there exists a sequence $w_{j}\to p$ with $\nu_{f(w_j)-s^*}\stackrel{w}\to\varepsilon_0$. Also $p_{i}\neq p,\quad,\forall i=1,2,\ldots$. From the definition of $T_f$, it now follows that $s^{*}\in T_{f}(p)$, i.e., $p_{k}\in T_{f}(p)$..  

\newtheorem{remarks}{Remarks}[section]
\begin{remarks}
Since the extremal points of $conT_{f}(p)$ are contained in $T_{f}(p)$ we conclude that $\psi_f$ is equal to the infimum of the probabilistic radius of $T_{f}(p)$ for $p$ in $V$, in other words:
$$\psi_{f}=\min_{p\in V}R_{T_{f}(p)}=\min_{p\in V}\min_{q\in T_{f}(p)}\nu_{f(p)-q}$$
\end{remarks}

The theorem below represents the probabilistic generalization of the famous Kakutani's fixed point theorem.

\begin{theo} 
Let $A$ be a nonempty, precompact, and convex subset of a PN-space $(V,\nu,\tau, \tau^*)$  and let $f\colon A\to A$  be any map. We set
$$\{(p,q)|p\in A,q\in T_{f}(p),\nu_{p-q}\stackrel{w}\to\varepsilon_0\}\subseteq V\times V$$
is closed. Then there is a point $p^{*}\in V$ such that $p^{*}\in T_{p}(p^*)$.
\end{theo}

Now, we consider our main theorem which shows that for any discontinuous function $f$ there exists at least an "approximate fixed point", i.e., a point $p^*$ whose image, $f(p^*)$, is no further away from $p^*$ than by a certain degree of measure of discontinuity for $f$.

\begin{theo}
Given a PN space $(V,\nu,\tau, \tau^*)$  . Let $A\subseteq V$   be a precompact and convex and $f\colon A\to A$ any map, then there exists a point $p^{*}\in A$ such that $\nu_{f(p^*)-p^*}\geq\psi_f$.
\end{theo}
{\bf Proof}\quad By Lemma 3.2, there exists $p^{*}\in conT_{f}(p^*)$. Hence 
$$\nu_{f(p^*)-p^*}(t)\geq \inf_{q \in T_{f}(p^*)}\nu_{f(p^*)-q}(t)\geq \inf_{p\in A}\inf_{q \in T_{f}(p)}\nu_{f(p)-q}(t)=\psi_{f}(t)$$
for every $t>0$, i.e., $\nu_{f(p^*)-p^*}\geq\psi_f$ as required.

\noindent{\bf Acknowledgement} 
The author would like to thank Prof. B. Schweizer and Prof. B. Lafuerza Guillen for sending their papers \cite{Alsina} and \cite{Guillen}.

\end{document}